\documentclass{ieeeaccess}
\usepackage{cite}

\usepackage{tikz}
\usetikzlibrary{math} 
\usetikzlibrary{calc}
\usetikzlibrary{arrows.meta, bending}
\NewSpotColorSpace{PANTONE}
\AddSpotColor{PANTONE} {PANTONE3015C} {PANTONE\SpotSpace 3015\SpotSpace C} {1 0.3 0 0.2}
\SetPageColorSpace{PANTONE}%
\def\centerarc[#1](#2)(#3:#4:#5);%
    {
    \draw[#1]([shift=(#3:#5)]#2) arc (#3:#4:#5);
    }
  
\usepackage{amsmath,amssymb,amsfonts}
\usepackage{graphicx}
\usepackage{textcomp}
\usepackage{booktabs}%
\usepackage{algorithm}%
\usepackage{algorithmicx}%
\usepackage{algpseudocode}%
\usepackage{listings}%

\usepackage{url}
\usepackage{makecell}
\usepackage{multirow}

\usepackage{siunitx}
\usepackage{xfrac}

\usepackage[hidelinks]{hyperref}
\usepackage[capitalise,nameinlink,sort]{cleveref}

\DeclareMathOperator*{\argmin}{arg\,min}

\usepackage[inline,shortlabels]{enumitem}

\usepackage[font={sf,scriptsize},           
            labelfont={bf,color=accessblue},
            caption=false]                  
           {subfig} 

\def\BibTeX{{\rm B\kern-.05em{\sc i\kern-.025em b}\kern-.08em
    T\kern-.1667em\lower.7ex\hbox{E}\kern-.125emX}}
\begin{document}
\history{Date of publication xxxx 00, 0000, date of current version xxxx 00, 0000.}
\doi{00.0000/ACCESS.0000.0000000}

\title{Swarm-Based Trajectory Generation and Optimization for Stress-Aligned 3D Printing}
\author{
\uppercase{Xavier Guidetti}\authorrefmark{1,2}, 
\uppercase{Efe C. Balta}\authorrefmark{1,2}, \IEEEmembership{Member, IEEE}, and 
\uppercase{John Lygeros}\authorrefmark{1}, \IEEEmembership{Fellow, IEEE}
}

\address[1]{Automatic Control Laboratory, ETH Zurich, Physikstrasse 3, 8092, Zurich, Switzerland (e-mail: \{xaguidetti,jlygeros\}@control.ee.ethz.ch)}
\address[2]{Inspire AG, Technoparkstrasse 1, 8005, Zurich, Switzerland (e-mail: efe.balta@inspire.ch)}

\tfootnote{Research supported by Innosuisse (project \textnumero 102.617 IP-ENG) and by the Swiss National Science Foundation under NCCR Automation (grant \textnumero 180545)}

\markboth
{Guidetti \headeretal: Preparation of Papers for IEEE TRANSACTIONS and JOURNALS}
{Guidetti \headeretal: Preparation of Papers for IEEE TRANSACTIONS and JOURNALS}

\corresp{Corresponding author: Efe C. Balta (e-mail: efe.balta@inspire.ch).}

\begin{abstract}
In this study, we present a novel swarm-based approach for generating optimized stress-aligned trajectories for 3D printing applications. 
The method utilizes swarming dynamics to simulate the motion of virtual agents along the stress produced in a loaded part. Agent trajectories are then used as print trajectories. With this approach, the complex global trajectory generation problem is subdivided into a set of sequential and computationally efficient quadratic programs. 
Through comprehensive evaluations in both simulation and experiments, we compare our method with state-of-the-art approaches. Our results highlight a remarkable improvement in computational efficiency, achieving a $115\times$ faster computation speed than existing methods. This efficiency, combined with the possibility to tune the trajectories spacing to match the deposition process constraints, makes the potential integration of our approach into existing 3D printing processes seamless. Additionally, the open-hole tensile specimen produced on a conventional fused filament fabrication set-up with our algorithm achieve a notable $\sim 10\%$ improvement in specific modulus compared to existing trajectory optimization methods.
\end{abstract}

\begin{keywords}
3D printing, additive manufacturing, fused filament fabrication, stress-aligned printing, swarming, trajectory optimization
\end{keywords}

\titlepgskip=-21pt

\maketitle

\section{Introduction}\label{sec:intro}


In 3D printing, a desired three-dimensional part is created by depositing material in a layer-wise fashion. This family of manufacturing processes, also known as Additive Manufacturing (AM), enables the creation of complex geometries and has been steadily gaining popularity in the last decades. The most widely AM technique is Fused Filament Fabrication (FFF) \cite{sculpteo2022}, which is also known as Fused Deposition Modeling (FDM) or Material Extrusion Additive Manufacturing \cite{ASTM_standard}. In FFF, a numerically controlled heated extruder moves along a predefined trajectory while depositing melted plastic, forming the desired part. Typically, the geometry to be manufactured is first sliced into equally spaced planes (known as \emph{planar} printing) \cite{gibson2021additive} or two-dimensional manifolds (known as \emph{non-planar} printing) \cite{nisja2021short, ahlers20193d}. Then, trajectories are created inside each slice. This process is generally achieved via a \emph{slicer} software, which mostly focuses on generating trajectories that form the part in a rapid or precise manner. Depending on the used material, desired geometry, printing equipment, and required geometrical accuracy, the slicing and trajectory generation process can be adapted to optimize the properties of the manufactured part \cite{hooshmand2021optimization,di2020reliable,delfs2016optimized,wang2020simultaneous}.

The recent developments of high performance materials, paired with a general technological improvement of the printing process, has enabled the usage of parts produced with FFF in a wide range of applications where high strength and stiffness are required \cite{brenken2018fused,jiang2017anisotropic}. Motivated by these developments, recent works have focused on generating stress-aligned printing trajectories that maximize the mechanical properties of FFF printed parts, both in the planar \cite{xia2020stress} and non-planar case \cite{fang2020reinforced, guidetti2022advanced}. Other works have tackled the same problem for the specific case of fiber-reinforced filaments, which must be printed in a continuous fashion \cite{barton2019fiber,chen2022field,yao20213d}. The general idea in this field of research is that, for a given material, manufacturing process, and geometry, there exists a stress-aligned trajectory that produces the desired part while improving its strength and stiffness under a general or specific load case. Finding such a trajectory is a challenging task due to complex part features, manufacturing process constraints and geometrical accuracy. The surveyed works that produce acceptable results in practice generally compute print trajectories by solving a very large optimization problem, which is complex and time-consuming. Additionally, these methods produce a single solution which cannot be tuned or refined to exactly match the printing process properties and mechanical requirements.

In this work, we propose a novel approach to the generation of stress-aligned print trajectories that is based on swarming dynamics. Each computed print line is the trace of a simulated agent moving through the part to be printed. Agents move in a swarm, and their motion is dictated by the stress generated by the part's load case, producing stress-aligned trajectories. By utilizing a swarm-based approach, the complex stress-aligned trajectory generation problem is broken down into a set of sequential and computationally efficient optimization problems. This significantly reduces the total time required for optimization. Additionally, the swarming dynamics enable great flexibility: our swarm-based trajectory generation method can be fine-tuned to produce trajectories that perfectly match the printing process limitations. The main contributions of this work are:
\begin{itemize}
    \item A novel, computationally efficient, and flexible swarm-based approach to stress-aligned trajectory generation for 3D printing, 
    \item A study of the methods' behavior with different settings, and
    \item The experimental comparison of swarm-based stress-aligned trajectory generation with a state-of-the-art method.
\end{itemize}

In \cref{sec:background} we introduce the methods of which we make use in our approach and discuss related literature. \cref{sec:method} details the swarm-based trajectory generation method. In \cref{sec:results} we analyze the performance of our method in detail and benchmark it. \cref{sec:concl} concludes the paper.

\section{Background}\label{sec:background}

\subsection{Finite Element Analysis} \label{sec:FEA}

Finite Element Analysis (FEA) is a very well known and widely used method to simulate numerically the stress field inside a loaded geometry \cite{reddy2019FEA}. We utilize it to compute the stress along which the printing trajectories for the optimized part will be aligned. In FEA, the part to be manufactured is first subdivided in a tetrahedral mesh. Then the forces produced by the predefined load case are added to the problem as boundary conditions. After assigning to the part its mechanical properties (i.e. Young's modulus and Poisson ratio), an FEA solver is used to find a numerical solution to the resulting set of differential equations. Solving this discretized problem corresponds to finding a Cauchy stress tensor $\boldsymbol\sigma$ for each node in the mesh. Typically, the Cauchy stress tensor is decomposed using eigenvalue decomposition to obtain the principal stresses and principal directions \cite{bower2009applied}. Now, the stress tensor can be represented as a diagonal matrix in a reference frame oriented along the principal directions of the decomposition. The diagonal entries $\sigma_{1}$, $\sigma_{2}$ and $\sigma_{3}$ of the principal stress matrix are ordered such that $|\sigma_{1}| \geq |\sigma_{2}| \geq |\sigma_{3}|$. In this work, we only consider the principal stress $\sigma_1$ and the corresponding normalized eigenvector $\mathbf{e}_1$. This eigenvector indicates the direction along which the principal stress acts. At any given node of the mesh, we indicate the product of the principal stress with the corresponding normalized eigenvector as $\mathbf{s} = \sigma_1\mathbf{e}_1$, and we loosely call it \emph{local principal stress vector}.


\subsection{Stress-Aligned Printing}


The goal of stress-aligned printing is to optimize the internal structure of a 3D printed part to enhance its mechanical performance. The approach exploits the anisotropic nature of AM processes. In particular, parts produced with FFF have been shown to be anisotropic, as their mechanical properties are better in the direction in which the plastic filament has been deposited. This effect can be moderately strong with materials such as PolyLactic Acid (PLA) (see Fig. 2 from \cite{fang2020reinforced}) or extremely strong with materials such as Liquid Crystal Polymers (LCP) \cite{gantenbein2018three}. In practice, after conducting a FEA for a given geometry and load case, the field of local principal stress vectors is used to guide the trajectory generation process. The objective is to produce trajectories for which the local printing direction is aligned with the local principal stress. This has been shown to maximize the strength and stiffness of the final object under the given load case \cite{gibson2021additive, fang2020reinforced, guidetti2023stress}. While the stress-alignment of trajectories has been the central goal of past works, it has been achieved at the expense of computational efficiency and of flexibility of the approach. Existing methods require long and complex computations, often relying on commercial optimization solvers, which makes them unsuitable for widespread use in conventional slicers. Furthermore, for a given geometry and load case, a unique trajectory is generated, with no possibility for refinement.
This is particularly limiting, since the manufacturability of trajectories is a central issue in FFF. Generally, a stress-aligned trajectory for a complex part is characterized by a variable spacing between the print lines. However, the possibility to deposit lines of variable width is constrained by the material properties and the extrusion process. Attempting to print beyond these constraints necessarily produces over- or under-extrusion, strongly reducing the quality of the printed part and its mechanical properties \cite{GuidettiIFAC,siqueira2017}. Thus, there exists a need for a computationally efficient stress-aligned trajectory generation algorithm that can be tuned to produce trajectories with a desired distribution of line spacing, ensuring manufacturability. This last aspect is becoming increasingly relevant as recent works have enhanced the quality and width range of variable width FFF \cite{guidetti2024data, guidetti2024force}.

\subsection{Swarming} \label{sec:swarming}

Numerous methods for swarm formation, modeling, and navigation for multiple autonomous agents have been proposed in the literature \cite{brambilla2013swarm, bayindir2016review, chung2018survey}. Often, the interactions between individuals (or between individuals and the space surrounding them) have been modeled using artificial potential functions. The main applications driving this line of research have been robot navigation and control \cite{khatib1986real,rimon1990exact,rimon1992} or multi-agent coordination \cite{reif1999social}. Interestingly, however, the inspiration for a more general approach to the study of swarm aggregation came from mathematical biology \cite{warburton1991tendency,grunbaum1994modelling}. Gazi and Passino \cite{gazi2004class,gazi2003stability,gazi2004stability} have introduced a class of attraction and repulsion functions between individuals that ensure the aggregation of a swarm moving through an environment, prevent individuals from making contact with each other, and allow for formation control. They have further extended their work to follow an energy approach \cite{gazi2013lagrangian}. This method, mimicking the behavior of swarms in nature, proposes that the motion of individual agents is dictated by a \emph{biological potential energy} to be minimized. The total energy of a swarm is given by the sum of its kinetic and potential energies. Trivially, a swarm kinetic energy corresponds to the sum of the kinetic energies of its individual components, as defined in classical mechanics. The notion of potential energy of a swarm, however, has been extended beyond classical physical potential energy to also include aggregation, environment and predator potentials. Given these notions, it is possible to model a swarm by applying a Lagrangian approach to the swarm total energy. Intuitively, the agents composing a swarm have a tendency to maintain their kinetic energy unchanged and to attain minimal potential. A flock of birds or a school of fish, for example, seek to travel with constant speed and direction (thus avoiding kinetic energy variations), to maintain a comfortable and safe distance between individuals (minimizing the aggregation potential), to explore areas rich in food (minimizing the environment potential), and to scatter and flee in the presence of threats (minimizing the predator potential).

\section{Method} \label{sec:method}

The approach we propose originates from the fascinatingly simple idea that the motion trajectories of a swarm of virtual point-mass agents can be used as manufacturing trajectories. If the swarm is carried through a mechanically loaded part by the load-induced stress flow, the resulting manufacturing trajectories are then stress-aligned. Similarly to other existing approaches in the literature \cite{fang2020reinforced, guidetti2023stress}, we separate the tasks of slicing the part and of generating trajectories on a slice. We assume that the slices are given and only consider the trajectory generation problem. Irrespective of whether slices are planar or not, the trajectories are thus generated over a 2D manifold.

The generation of trajectories for FFF has specific requirements and peculiarities that make the existing approaches to swarm modeling not immediately usable for the task. The method introduced in \cref{sec:swarming} must be modified or extended in three areas: 
\begin{enumerate}
    \item agents in the swarm must avoid moving through existing trajectories, \label{enum:crossing}
    \item the set of trajectories must cover the entire part, and \label{enum:coverage}
    \item agents can be added to or removed from the swarm at will. \label{enum:spawning}
\end{enumerate}

Condition \ref{enum:crossing} is not enforced in classical swarm modeling, since physical agents must simply avoid collisions among each other. Only the present location of agents is relevant for planning, and agents are allowed to cross existing trajectories. In swarming for manufacturing, however, both a spatial and temporal separation must be maintained between agents: the manufacturing trajectories cannot overlap as depositing material twice at the same location creates defects in the part.

Condition \ref{enum:coverage} ensures that the outside of the part produced with the generated trajectories matches the desired geometry, and that no voids are left inside the part. Additionally, in FFF (and in AM in general) the rate of material deposition is bounded. For this reason, good quality coverage of the part is achieved when the spacing between neighboring manufacturing trajectories respects upper and lower bounds. If trajectories are too far apart the space in between cannot be entirely filled, if trajectories are too close together too much material will be deposited at the same location.

Finally, Condition \ref{enum:spawning} exploits the non-physical nature of agents in the manufacturing problem to simplify the part coverage problem (Condition \ref{enum:coverage}). Contrarily to a swarm of physical agents, we can change the number of agents in the swarm at any time. For example, it is possible to \emph{spawn} new agents (which will produce a new manufacturing trajectory) where other trajectories have diverged, or to \emph{kill} existing agents (which will interrupt an existing manufacturing trajectory) where trajectories compress excessively.

\subsection{Swarm-Based Trajectory Planning for Manufacturing}

Our approach has been principally inspired by the notions of aggregation and environment potentials from \cite{gazi2013lagrangian}, which we have adapted to the context of manufacturing. Let us consider a mechanically loaded part for which a FEA has been conducted as discussed in \cref{sec:FEA}. We create a swarm of agents at the location where the largest external force is applied to the part. At initialization, the agents are homogeneously spaced, and neighboring agents are at a predefined distance matching the nominal desired distance between FFF print lines. We define the swarm potential as
\begin{equation}
    P = P_a + K P_e \,, \label{eq:swarm_pot}
\end{equation}
where $P_a$ and $P_e$ are the aggregation and environment potentials of the swarm, and $K$ is a tuning constant. $P_a$ encodes the quality of the swarm formation: it is minimal when the spacing between neighboring agents corresponds to a predefined desired distance, and increases in case of spacing deviations. $P_e$ quantifies the stress alignment of the agents' motion: it is minimal when agents follow the stress flow exactly, and increases when they deviate from it.
To have the swarm draw a well-spaced set of stress-aligned trajectories through the part, we propose to utilize \cref{alg:outline}.

\begin{algorithm}
\caption{Swarm-Based Trajectory Generation}\label{alg:outline}
\begin{algorithmic}[1]
\Require Agents initial location, Step size, Agents desired distance, FEA simulation, $K$
\While{Part is not fully covered}
\State Advance each agent by one step size along the local principal stress direction
\State Reposition agents by minimizing $P$
\If{Two agents are too close together}
        \State Kill one of them
\EndIf
\If{Two agents are too far apart}
        \State Spawn one agent in between
\EndIf
\EndWhile

\end{algorithmic}
\end{algorithm}

\subsection{Methodological Details} \label{sec:implementation}

\cref{alg:outline} is obviously only an outline of the method, that we first introduce in a simplified manner for clarity. In this section, we describe every required detail that composes our method for swarm-based trajectory planning for manufacturing.

\subsubsection{Step Direction and Size}

The solution to an FEA simulation of stresses in a mechanical component produces a stress flow in the part. As the principal stress vectors forming the stress flow originate from the eigenvectors of the Cauchy stress tensor, only their direction is relevant to the stress alignment problem, while their orientation is not. Intuitively, local alignment between a print line and a stress vector is left unchanged by flipping the vector orientation. To produce a stress flow free from this \emph{\SI{180}{\degree} ambiguity} (i.e. heterogeneous stress vectors orientations), past approaches have used simulated annealing \cite{metcalf1994resolving} and rectification along the main Cartesian component \cite{guidetti2023stress}. In this work, we exploit the particle swarming nature of the approach to solve the orientation ambiguity by using the momentum of the agents. At each iteration, an individual agent is displaced according to the local principal stress vector \emph{direction}. However, we select the displacement \emph{orientation} that best aligns with the displacement of the agent in the previous iteration. In practice, using this approach, no sharp changes (i.e. larger than \SI{90}{\degree}) in an agent's trajectory are possible.

Contrarily to intuition, the size of the agents steps does not depend on the magnitude of the stress vectors. To ensure that agents advance in a uniform front, every agent moves by a predefined and fixed step size. The uniformity of the front is a condition required to simplify the optimization problem, as we will see in the rest of this section. To ensure that no trajectories overlap (see \cref{sec:method}, Condition \ref{enum:crossing}), we set the step size $h = l$, where $l$ is the desired distance between neighboring agents.

To summarize, we consider an agent $x_1$ that at iteration $k$ is located in $x_1(k)$ and the local principal stress vector $\mathbf{s}_1(k)$ (which we project on the slice if necessary). We indicate with $\hat{\mathbf{s}}_1(k) = \sfrac{\mathbf{s}_1(k)}{\|\mathbf{s}_1(k)\|_2}$ the normalized local principal stress vector. After one step, the agent will be located in
\begin{equation}
	t_1(k+1) = x_1(k) + \mu \hat{\mathbf{s}}_1(k) h \,, \label{eq:step}
\end{equation}
where $\mu$ is the momentum term used to circumvent the ambiguity problem: we set $\mu=1$ when $\mathbf{s}_1(k) \cdot (x_1(k)-x_1(k-1)) \geq 0$ and $\mu=-1$ otherwise. $x_1(k-1)$ naturally denotes the location of the agent at the previous iteration of the algorithm, and $x_1(k)-x_1(k-1)$ corresponds to its last displacement. The forward step motion of two agents can be observed in \cref{fig:P_a}.

\subsubsection{Agents Potential Functions} \label{sec:potential}


Two potential functions are required to obtain the swarm potential introduced in \cref{eq:swarm_pot}: the environment potential $P_e$, and the aggregation potential $P_a$. The potential functions commonly utilized in the literature \cite{gazi2013lagrangian} are nonlinear, and the resulting total potential of a swarm is generally non-convex. Optimizing such functions remains a complex, inefficient, and time-consuming task \cite{danilova2022recent, snyman2005practical}. To simplify and accelerate the potential optimization problem, we introduce two quadratic environment and aggregation potential functions. Using these, \cref{eq:swarm_pot} becomes a quadratic programming (QP) problem \cite{frank1956algorithm}, which can be solved very efficiently by numerous existing solvers.

\begin{figure}[htbp]
\centering
\begin{tikzpicture}
    
    \draw[thick,teal,-{Stealth[length=3mm]}] (3.8,0) -- (3.5,2) node[at start,black,anchor=east]{$x_1(k-1)$} node[anchor=east,black]{$x_1(k)$};
    \draw[thick,teal,-{Stealth[length=3mm]}] (6,0) -- (5.8,2.2) node[at start,black,anchor=west]{$x_2(k-1)$} node[anchor=west,black]{$x_2(k)$};

    \draw[thick,black,-{Stealth[length=3mm]}] (3.5,2) -- (1.5, 5) node[anchor=east,black]{$t_1(k+1)$}
    node[midway,below,sloped]{$\hat{\mathbf{s}}_1(k)h$};
    \draw[thick,black,-{Stealth[length=3mm]}] (5.8,2.2) -- (6.4, 6) node[anchor=west,black]{$t_2(k+1)$}
    node[midway,below,sloped]{$\hat{\mathbf{s}}_2(k)h$};

    \draw[thick,orange,-{Stealth[length=3mm]}] (1.5, 5) -- (2,6.1) node[anchor=south east,black]{$x_1(k+1)$} node[midway, above,sloped] {$\delta_1$};
    \draw[thick,orange,-{Stealth[length=3mm]}] (6.4, 6) -- (6, 7.5) node[anchor=west,black]{$x_2(k+1)$} node[midway, above,sloped] {$\delta_2$};

    \draw[thick,purple,-{Stealth[length=3mm]}] (2,6.1) -- (6, 7.5)
    node[midway,above,sloped]{$\mathbf{v}$};
    \tikzmath{\x1 = 3.4; \y1 =0.8; \x2 = \x1 + 0.65*4.2; \y2 =\y1 + 0.65*0.5;} 
    \draw[thick,blue,-{Stealth[length=3mm]}] (\x1,\y1) -- (\x2, \y2)
    node[midway,above,sloped]{$\hat{\mathbf{d}}$};

    \draw[black, very thin, fill=gray!5] (0.5,-3) rectangle (7.5,-0.7);
    \draw[thick,-{Stealth[length=3mm]}] (2+0.5,-3+0.5) -- (2+0.5+4.12,-3+0.5+0.4905) node[below]{$\textit{proj}_{\hat{\mathbf{d}}}\mathbf{v}$};
    \draw[thick,-{Stealth[length=3mm]}] (2+0.5,-3+0.5) -- (6+0.5-4.12,-1.5+0.5-0.4905) node[near end,left]{$\textit{oproj}_{\hat{\mathbf{d}}}\mathbf{v}$};
    
    \draw[thick,purple,-{Stealth[length=3mm]}] (2+0.5,-3+0.5) -- (6+0.5,-1.5+0.5) node[midway,above,sloped]{$\mathbf{v}$};
    \draw[thick,blue,-{Stealth[length=3mm]}] (2+0.5,-3+0.5) -- (\x2-\x1+2+0.5,\y2-\y1-3+0.5) node[midway,below,sloped]{$\hat{\mathbf{d}}$};

    \draw[dashed] (2+0.5+4.12,-3+0.5+0.4905) -- (6+0.5,-1.5+0.5);
    \draw[dashed] (6+0.5-4.12,-1.5+0.5-0.4905) -- (6+0.5,-1.5+0.5);

\end{tikzpicture}
\caption{Notation of the agents locations and distance vectors utilized in the definition of the environment potential $P_e$ and of the aggregation potential $P_a$}
\label{fig:P_a}
\end{figure}
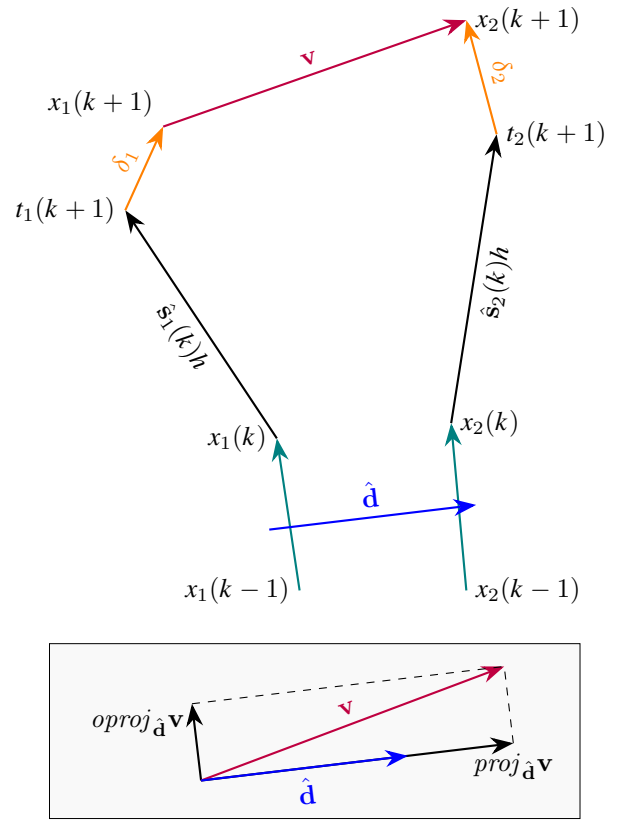

To define the environment potential $P_e$, let us consider an individual agent $x_1$ that has  advanced by a step in the local stress direction following \cref{eq:step} and is now located in $t_1(k+1)$. We can imagine that $t_1(k+1)$ is the ``ideal'' location for the agent, as the trajectory of the agent between $x_1(k)$ and $t_1(k+1)$ is perfectly aligned with the local stress vector $\mathbf{s}_1(k)$. We assign the agent a \emph{virtual mass} $m_1(k+1) \propto \|\mathbf{s}_1(k)\|_2$, that is proportional to the local stress vector magnitude\footnote{In practice, we use the largest principal stress vector found in the FEA simulation as a normalization factor. We set, for each agent $x_i(k+1)$ and corresponding local stress vector $\mathbf{s}_i(k)$, a virtual mass $m_i(k+1) = \sfrac{\|\mathbf{s}_i(k)\|_2}{\max\|\mathbf{s}\|_2}$}. We can now imagine the environment potential associated to the agent as the \emph{energy} required to reposition the agent away from its ideal location $t_1(k+1)$, and bring it to a final location $x_1(k+1)$, as visible in \cref{fig:P_a}. We define this as
\begin{equation}
	P_e\bigl(x_1(k+1)\bigr) =  m_1(k+1) \|x_1(k+1)-t_1(k+1)\|_2 \,, \label{eq:P_e}
\end{equation}
which is a quadratic function with a global minimum $P_e\bigl(x_1(k+1)\bigr)=0$ when the agent is not repositioned. Using the virtual mass term $m$, we make the repositioning of agents with a larger associated stress more costly. This ensures that where the local stress is larger, the alignment between trajectories and stress will be higher.

The aggregation potential $P_a$ is computed by comparing the locations of neighboring agents. As discussed in \cref{sec:method}, Condition \ref{enum:crossing}, we want to avoid agent trajectories that overlap or cross. An intuitive way to achieve this would be to encode existing trajectories as locations to be avoided by the moving agents (for example, by including them in the environment potential). However, this approach presents a major scalability issue. As the number of iteration of \cref{alg:outline} grows, the size of existing trajectories keeps on increasing, making the size and complexity of the optimization problem larger and larger, and leading to computation issues. Additionally, this problem is amplified if the parts to be manufactured are large with respect to agent spacing. We avoid this issue altogether by ensuring that the moving agents never encounter any past trajectory. In this case, their potential functions can be designed to only depend on the location of other agents at the \emph{current} iteration, and not at \emph{all past} iterations. Instead of evaluating the distance between agents in a straight line, we decompose it into \emph{radial distance}, which is measured perpendicularly to the displacement direction of the agents, and \emph{axial distance}, measured in the displacement direction of the agents.
We achieve an orderly advance of the agents by keeping neighboring agents at a desired radial distance and by minimizing their axial distance. Thus, the swarm advances as a front of homogeneously spaced agents, which are aligned along a quasi-straight front line. 
We consider two neighboring agents $x_1$ and $x_2$ whose locations at iteration $k+1$ are denoted as $x_1(k+1)$ and $x_2(k+1)$. We define the vectors
\begin{equation}
    \mathbf{v} = x_2(k+1) - x_1(k+1) \,,
\end{equation}
corresponding to the directed distance between the agents, and
\begin{align}
    \mathbf{d} &= -\bigl( x_1(k) - x_1(k-1) + x_2(k) - x_2(k-1) \bigr) ^\perp \,, \label{eq:d} \\
    \hat{\mathbf{d}} &= \frac{\mathbf{d}}{\|\mathbf{d}\|_2} \, ,
\end{align} 
the unit vector orthogonal\footnote{The notation $\cdot^\perp$ used in \cref{eq:d} corresponds, in two dimensions, to a counterclockwise rotation of the vector by \SI{90}{\degree}, i.e. $\mathbf{a}^\perp = \big[\begin{smallmatrix}
  0 & -1\\
  1 & 0
\end{smallmatrix}\big] \mathbf{a}$.} to the averaged direction of travel of the two agents in the previous iteration. Both vectors are shown in \cref{fig:P_a} for clarity.
We define the aggregation potential $P_a$ of agent $x_1$ with respect to its neighbor $x_2$ as
\begin{equation}
	P_a\bigl(x_1(k+1)|_{x_2}\bigr) = (\|\textit{proj}_{\hat{\mathbf{d}}}\mathbf{v}\|_2 - l)^2 + \|\textit{oproj}_{\hat{\mathbf{d}}}\mathbf{v}\|_2^2 \,, \label{eq:P_a}
\end{equation}
where the norm in the first term, encoding the requirement to keep the radial distance close to the desired distance $l$, is computed as
\begin{equation}
	\|\textit{proj}_{\hat{\mathbf{d}}}\mathbf{v}\|_2 = \mathbf{v}\cdot\hat{\mathbf{d}} \,, \label{eq:proj}
\end{equation}
and the norm in the second term, used to minimize the axial distance, is computed as
\begin{equation}
	\|\textit{oproj}_{\hat{\mathbf{d}}}\mathbf{v}\|_2 = \|\mathbf{v} - (\mathbf{v}\cdot\hat{\mathbf{d}})\hat{\mathbf{d}}\|_2 \,. \label{eq:oproj}
\end{equation}
As $\hat{\mathbf{d}}$ and $l$ are fixed, $P_a$ is a quadratic function depending on $x_1(k+1)$ and $x_2(k+1)$, and has a global minimum $P_a\bigl(x_1(k+1)|_{x_2}\bigr) = 0$ when the two agents advance side by side at distance $l$.

\subsubsection{Agents Adjacency} \label{sec:adjacency}

The aggregation potential $P_a$ of the swarm is obtained by comparing the locations of different agents in the swarm. We greatly simplify the computation task by only comparing \emph{neighboring} agents. Every agent has two neighbors, each being the closest agent in either orientation along the \emph{radial} direction. In simpler terms, agent $x_i$ has neighbors $x_{i-1}$ and $x_{i+1}$; any other agent in the swarm has no aggregation potential with respect to $x_i$. The adjacency of agents is fixed since initialization and does not change during their motion through the part. Given the commutative nature of $P_a$ as defined in \cref{eq:P_a}, we consider as adjacent during the computations only agents $x_i$ and $x_{i+1}$. This corresponds to having $P_a\bigl(x_i|_{x_{j\neq i+1}}\bigr) = 0$.

\subsubsection{Boundary Conditions}

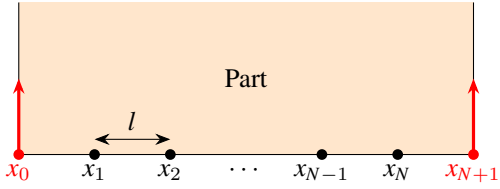
\begin{figure}[htbp]
\centering
\begin{tikzpicture}
    \fill[orange!20] (0,0) rectangle (6,2);
    \draw (0,0) -- (6,0);
    \draw (0,0) -- (0,2);
    \draw (6,0) -- (6,2);
    \node at (3,1) {Part};
    
    \fill[red] (0,0) circle (2pt) node[below] {$x_0$};
    \fill (1,0) circle (2pt) node[below] {$x_1$};
    \fill (2,0) circle (2pt) node[below] {$x_2$};
    \fill (3,0) circle (0pt) node[below] {$\cdots$};
    \fill (4,0) circle (2pt) node[below] {$x_{N-1}$};
    \fill (5,0) circle (2pt) node[below] {$x_N$};
    \fill[red] (6,0) circle (2pt) node[below] {$x_{N+1}$};
    
    \draw[{Stealth[length=2mm]}-{Stealth[length=2mm]}] (1,0.2) -- (2,0.2) node[midway, above] {$l$};
    \draw[-{Stealth[length=2mm]},red, very thick] (0,0) -- (0,1);
    \draw[-{Stealth[length=2mm]},red,very thick] (6,0) -- (6,1);
    
\end{tikzpicture}
\caption{Initial swarm of equally spaced agents $x_{1,\ldots,N}$ and boundary agents $x_0$ and $x_{N+1}$. The displacement of the boundary agents is one dimensional, as they are constrained to follow the part boundary.}
\label{fig:bound}
\end{figure}

To produce parts with a desired shape, the generated trajectories must cover the part entirely, as discussed in \cref{sec:method}, Condition \ref{enum:coverage}. By simply defining a swarm of initial agents and making it propagate through the part as explained in \cref{alg:outline}, complete coverage is not guaranteed. In particular, without further constraints, agents are free to step outside the part or to leave an empty non-covered space between the part boundaries and the swarm. To avoid such scenarios, we introduce \emph{boundary conditions} to constrain and condition the swarm.
As the swarm is initialized, at the beginning of \cref{alg:outline}, two \emph{boundary agents} are added at the extremities of the swarm. In a swarm $\mathcal{X}$ composed of $N$ agents $x_{1,\ldots,N}$, we create the boundary agents $x_0$ and $x_{N+1}$, as shown in \cref{fig:bound}. The motion of these two agents is constrained to follow the part boundary. As a consequence, when computing the swarm trajectories, the agents $x_1$ and $x_N$ maintain approximately a distance $l$ from the part boundaries, ensuring complete coverage and avoiding boundaries crossing.

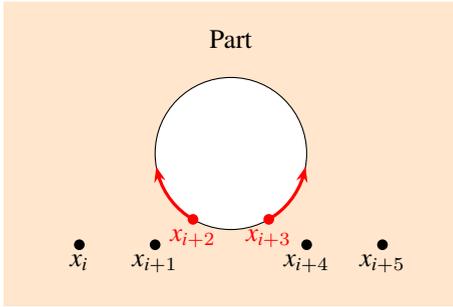
\begin{figure}[htbp]
\centering
\begin{tikzpicture}
    \fill[orange!20] (0,0) rectangle (6,4);
    \fill[white] (3,2) circle (1);
    \draw[black, ] (3,2) circle (1);
    \node at (3,3.5) {Part};
    
    \def\angle{30}
    \def\length{1}
    
    \fill (1,0.8) circle (2pt) node[below] {$x_i$};
    \fill (2,0.8) circle (2pt) node[below] {$x_{i+1}$};
    \fill (4,0.8) circle (2pt) node[below] {$x_{i+4}$};
    \fill (5,0.8) circle (2pt) node[below] {$x_{i+5}$};

    \coordinate (start1) at (3-sin{\angle},2-cos{\angle});
    \coordinate (end1) at ($(start1) + (180-\angle:\length)$);
    \coordinate (start2) at (3+sin{\angle},2-cos{\angle});
    \coordinate (end2) at ($(start2) + (\angle:\length)$);
    \fill[red] (start1) circle (2pt) node[below] {$x_{i+2}$};
    \fill[red] (start2) circle (2pt) node[below] {$x_{i+3}$};    
    \centerarc[-{Stealth[length=2mm]},red,very thick](3,2)(\angle-90:\angle-40:1);
    \centerarc[{Stealth[length=2mm]}-,red,very thick](3,2)(-\angle-90-50:-\angle-90:1);

\end{tikzpicture}
\caption{Internal boundary agents $x_{i+2}$ and $x_{i+3}$ added to the swarm when encountering a hole in the part}
\label{fig:inner_bound}
\end{figure}

During trajectory generation, the swarm of agents can encounter other part boundaries (such as holes or other non-convex features). In this case, two additional \emph{internal} boundary agents are added to the swarm in between the two agents closest to the boundary, as shown in \cref{fig:inner_bound} The agents adjacency is updated accordingly, and the trajectory generation according to \cref{alg:outline} continues. Should the two internal boundary agents meet along their motion (for example, when the swarm has moved past a hole), they are removed from the swarm.

\subsubsection{Optimization Problem}

To compute the total potential of the swarm $P$, that we intend to minimize at each iteration, we first need to define the total environment and aggregation potentials of the swarm. Considering a swarm $\mathcal{X}$ of $N$ agents $x_{1,\ldots,N}$, at iteration $k+1$, the swarm environment potential is simply
\begin{equation}
    P_e\bigl(\mathcal{X}(k+1)\bigr) = \sum_{i=1}^N P_e\bigl(x_i(k+1)\bigr) \,,
\end{equation}
where the environment potential of an individual agent is given in \cref{eq:P_e}. Similarly, we define the swarm aggregation potential as
\begin{equation}
    P_a\bigl(\mathcal{X}(k+1)\bigr) = \sum_{i=0}^N P_a\bigl(x_i(k+1)|_{x_{i+1}}\bigr) \,,
\end{equation}
where the aggregation potential of an individual agent is given in \cref{eq:P_a} and includes the notion of adjacency introduced in \cref{sec:adjacency}.

\begin{figure}[htbp]
\centering
\begin{tikzpicture}
    \def\sx{5.5}
    \def\sy{1}
    
    \draw[magenta, very thin, fill = magenta!20,rotate around={atan(\sy/\sx):(\sx,\sy)}] (3/4*\sx,-1/8*\sx+\sy) rectangle ++(\sx/2,\sx/4);
    \draw[{Stealth[length=2mm]}-{Stealth[length=2mm]},rotate around={atan(\sy/\sx):(\sx,\sy)}] (5/4*\sx+0.2,-1/8*\sx+\sy) -- ++(0,\sx/4) node[midway, right]{$h/4$};
    \draw[{Stealth[length=2mm]}-{Stealth[length=2mm]},rotate around={atan(\sy/\sx):(\sx,\sy)}] (3/4*\sx,-1/8*\sx+\sy-0.2) -- ++(\sx/2,0) node[midway, below]{$h/2$};

    \draw[thick,black,-{Stealth[length=3mm]}] (0,0) -- (\sx, \sy) node[anchor=north,black]{$t_1(k+1)$} node[at start,black,anchor=south]{$x_1(k)$}
    node[midway,above,sloped]{$\hat{\mathbf{s}}_1(k)h$};

    \draw[thick,orange,-{Stealth[length=3mm]}] (\sx, \sy) -- ++(0.4,0.7) node[anchor=south,black]{$x_1(k+1)$} node[near start, above,sloped] {$\delta_1$};
\end{tikzpicture}
\caption{Representation of the constraints used in the optimization problem \eqref{eq:opt_problem} to limit the displacement of the agents around their ideal location $t_i(k+1)$. The two constraints form the box depicted in magenta, which constrains $x_1(k+1)$.}
\label{fig:const}
\end{figure}
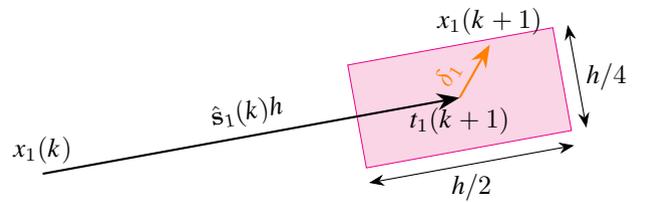

For a given tuning hyperparameter $K$, we can formulate the optimization problem being solved at each iteration of \cref{alg:outline} as
\begin{align}
\min_{x \in \mathcal{X}}\quad &P_a\bigl(\mathcal{X}(k+1)\bigr) + K  P_e\bigl(\mathcal{X}(k+1)\bigr) \label{eq:opt_problem} \\ 
\text{s.t.} \quad & \|\textit{proj}_{\hat{\mathbf{s}}_1}\delta_i(k+1)\|_2 \leq h/4  \qquad i = 1,\ldots,N \nonumber \\
& \|\textit{oproj}_{\hat{\mathbf{s}}_1}\delta_i(k+1)\|_2 \leq h/8  \qquad i = 1,\ldots,N \nonumber \\
& x_i(k+1) \in \partial P \qquad \text{for boundary agents} \,,\nonumber 
\end{align}
where $\delta_i(k+1) = \bigl(x_1(k+1)-t_i(k+1)\bigr)$. The definitions of the projection terms in the constraints correspond to those given in \cref{eq:proj} and \cref{eq:oproj}, and $\partial P$ indicates the boundary of the part. The inequality constraints limit the displacement of agents to a box around their ideal location $t_i(k+1)$, as shown in \cref{fig:const}, and are necessary to ensure that agents keep advancing and do not cross. With a suitable change of coordinates, all constraints can be transformed into input bounds in the form $a_i \leq x_i(k+1) \leq b_i$. The resulting problem is a bound-constrained QP, which can be solved very efficiently.

\subsubsection{Killing and Spawning} \label{sec:spawnkill}

During trajectory generation, we exploit the fact that agents can be spawned or killed when necessary (see \cref{sec:method}, Condition \ref{enum:spawning}). The procedure we utilize is based on the swarm potential $P$, following again the intuition that the ideal swarm behavior is achieved by minimizing $P$. At every iteration of \cref{alg:outline}, we study the swarm (which we denote in this section with $\mathcal{X}_0(k+1)$) to identify the agent with the closest neighbor and the agent with the furthest neighbor. This corresponds to locating the regions where the trajectories are most compressed and most expanded. In the first case, we kill the \emph{compressed} agent to form the swarm $\mathcal{X}_{-1}(k+1)$; in the second case, we spawn a new agent in the empty space between neighbors to form the swarm $\mathcal{X}_{+1}(k+1)$. We then solve the optimization problem \eqref{eq:opt_problem} for $\mathcal{X}_0(k+1)$, $\mathcal{X}_{-1}(k+1)$, and $\mathcal{X}_{+1}(k+1)$. We denote as $\mathcal{P}\bigl(\mathcal{X}(k+1)\bigr)$ the solution of \cref{eq:opt_problem} for a swarm $\mathcal{X}(k+1)$. Finally, after computing the optimized potential of each swarm, and normalizing it with the number of agents in each swarm, we can make a comparison. We select and keep for the next iteration the swarm 
\begin{equation}
    \mathcal{X}(k+1) = \argmin_{q=-1,0,1} \frac{\mathcal{P}\bigl(\mathcal{X}_q(k+1)\bigr)}{|\mathcal{X}_q|} \,,
\end{equation}
which has the lowest potential among the evaluated scenarios. Clearly, this approach relies on heuristics which make only one spawning or one killing possible at each iteration. This choice was done to reduce the computational burden, as it allows \cref{alg:outline} to change the number of agents in the swarm in a principled way, while only solving \cref{eq:opt_problem} three times per iteration. 


\section{Results} \label{sec:results}

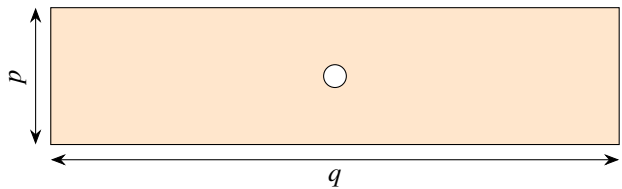
\begin{figure}[htbp]
\centering
\begin{tikzpicture}
    \def\tlen{7.5}
    \fill[orange!20] (0,0) rectangle (\tlen,0.24*\tlen);
    \draw[black] (0,0) rectangle (\tlen,0.24*\tlen);
    \fill[white] (\tlen/2,0.24/2*\tlen) circle (0.24*\tlen/6/2);
    \draw[black] (\tlen/2,0.24/2*\tlen) circle (0.24*\tlen/6/2);

    \draw[{Stealth[length=2mm]}-{Stealth[length=2mm]}] (-0.2,0) -- (-0.2,\tlen*0.24) node[midway, sloped, above] {$p$};
    \draw[{Stealth[length=2mm]}-{Stealth[length=2mm]}] (0,-0.2) -- (\tlen,-0.2) node[midway, below] {$q$};
    
\end{tikzpicture}
\caption{Open-hole tensile specimen according to \cite{american2023standard}. The dimensions are $p = \SI{36}{mm}$ and $q = \SI{150}{mm}$, and the thickness is \SI{2}{mm}. The hole has a diameter of \SI{6}{mm}. The specimen is evaluated in a tensile strength test, with tension applied at the two narrow extremities.}
\label{fig:specimen}
\end{figure}

In this section, we test and benchmark the swarm-based trajectory generation algorithm on a demonstrator part. The part we selected is a loaded open-hole tensile specimen. The specimen is taken from the ASTM standard for open-hole tensile strength of polymer matrix composite laminates \cite{american2023standard}, and described in \cref{fig:specimen}. We first analyze the behavior of \cref{alg:outline} by studying the trajectories generated for the specimen. Then, we conduct an experiment to quantify the improvements obtained in the specimen strength via the proposed method.

\subsection{Trajectories Generation}

We implement the swarm-based trajectory generation algorithm in Python (including SciPy \cite{2020SciPy-NMeth} and NumPy \cite{numpy}); the optimization problems are solved efficiently by utilizing CasADi \cite{Andersson2018} as an interface and OSQP \cite{osqp} as a solver. We first discuss the effect of the hyperparameter $K$, and then compare the trajectory generation performance with the method from \cite{guidetti2023stress}.

\subsubsection{Effect of $K$}
\begin{figure}
\centering
\subfloat[$K=0.5$\label{fig:traj05}]{\includegraphics[trim={6mm 6mm 6mm 6mm},clip]{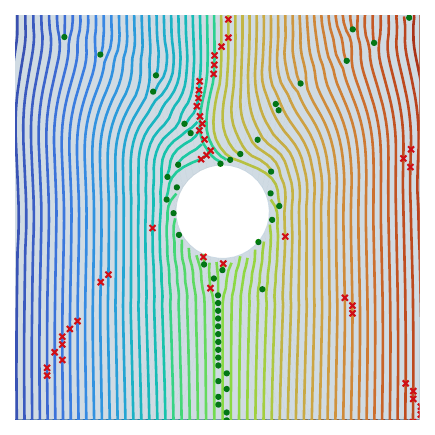}}
\vfill
\subfloat[$K=5$\label{fig:traj5}]{\includegraphics[trim={6mm 6mm 6mm 6mm},clip]{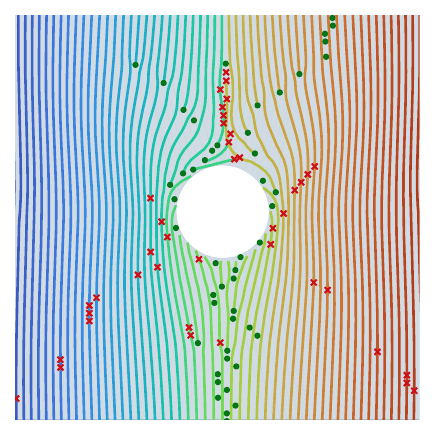}}
\vfill
\subfloat[$K=50$\label{fig:traj50}]{\includegraphics[trim={6mm 6mm 6mm 6mm},clip]{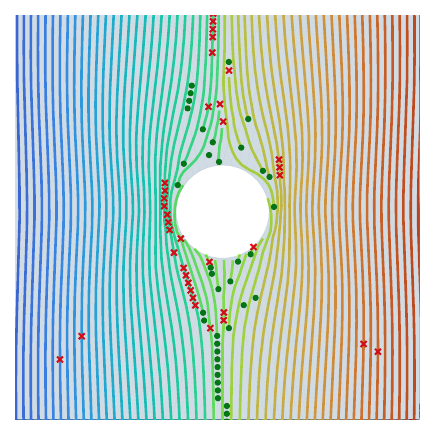}}
\caption{Comparison of manufacturing trajectories for the open-hole tensile specimen generated with different values of the hyperparameter $K$. The swarm advances from the bottom of the figures towards the top. Green circles and red crosses indicate the locations in which \cref{alg:outline} evaluates a potential spawn or kill, respectively (see \cref{sec:spawnkill}).}
\label{fig:K_comp}
\end{figure}

In \cref{eq:opt_problem}, $K$ modifies the relative weight of the aggregation potential $P_a$ with respect to the environment potential $P_e$. Selecting smaller values of $K$ in the optimization produces swarms that mostly minimize $P_e$: this corresponds to better swarming behavior (i.e. agents advancing as a front and evenly spaced) at the expense of stress alignment. Conversely, larger values of $K$ increase the importance of $P_a$ in the optimization problem: the agents follow the stress in the part more, but the swarming behavior is worsened. In \cref{fig:K_comp} we show the trajectories generated around the hole of the tensile specimen for three different values of $K$. With $K=0.5$ (\cref{fig:traj05}) the trajectories are uniformly spaced, but they follow the stress circulating around the hole poorly, which reduces the manufactured part strength. With $K=50$ (\cref{fig:traj50}) the trajectories follow the stress around the hole, but become very compressed in the regions to its left and right, which can cause defects in the manufacturing process. The intermediate case, produced with $K=5$ (\cref{fig:traj5}) appears to constitute a satisfactory trade-off, with the trajectories showing a better stress following behavior than in the $K=0.5$ case, and a better swarming behavior than in the $K=50$ case.

\begin{figure}[htbp]
\centering
\includegraphics{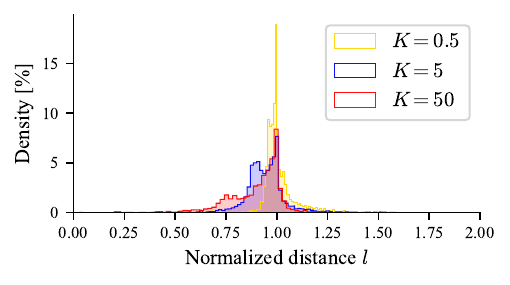}
\caption{Distribution of distances between trajectories for different values of $K$. The distance has been normalized using the nominal line spacing $l$.}
\label{fig:distrib}
\end{figure}

To quantify the effect of $K$, we analyze stress alignment and trajectory spacing. For stress alignment, we introduce the metric
\begin{equation}
    \bar{\beta} = \frac{\sum_{z=1}^Z m_z \|\hat{\mathbf{s}}_z \cdot \mathbf{p}_z\|}{\sum_{z=1}^Z m_z} \,,
\end{equation}
where $z$ indexes all the $Z$ points in a set of generated trajectories, and $\hat{\mathbf{s}}_z$ and $\mathbf{p}_z$ are the normalized local stress and the printing direction at the point $z$. The alignment metric $\bar{\beta}$ is weighted by the magnitude of the local stress (encoded in the virtual mass $m_z$, see \cref{sec:potential}), and  can range between $0$ (no alignment) and $1$ (perfect alignment). The alignment values for the three studied cases are given in \cref{tab:K}. For trajectory spacing, we compute the distance of every point in the trajectory set from the following trajectory, and visualize the results as a distribution in \cref{fig:distrib}. We also compute the variance of each distribution and report it in \cref{tab:K}. Both metrics are in agreement with the qualitative intuitions obtained from \cref{fig:K_comp}, as both $\bar{\beta}$ and the variance of the distance distribution increase monotonically with $K$.

\begin{table*}[htbp]
\centering
\caption{Stress alignment and variance of the distribution of distances of the swarm-based trajectory generation method (for different values of $K$) and of the global trajectory generation method from \cite{guidetti2023stress}}
\renewcommand{\arraystretch}{1.3}
\renewcommand\cellset{\renewcommand\arraystretch{0.8}%
\setlength\extrarowheight{0pt}}
\begin{tabular}[h]{@{}l c c c c c@{}}
\toprule
& \multicolumn{3}{c}{Swarm-based method} && \multirow{ 2}{*}{\makecell[c]{Global method\\from \cite{guidetti2023stress}}} \\
\cmidrule{2-4}
& $K=0.5$ & $K=5$ & $K=50$ \\
\midrule
$\bar{\beta}$ & \num{0.981} & \num{0.993} & \num{0.998} && \num{0.983} \\
Variance & \SI{6.1e-3}{} & \SI{12.9e-3}{} & \SI{16.4e-3}{} && \SI{4.4e-4}{} \\
\bottomrule
\end{tabular}
\label{tab:K}
\end{table*}

\subsubsection{Comparison and Benchmarking}

We compare \cref{alg:outline} with the method proposed in \cite{guidetti2023stress}, a state-of-the-art method for stress-aligned trajectory generation for 3D printing, which we call \emph{global trajectory generation}. We utilize both algorithms to generate print trajectories for one slice of the open-hole tensile specimen detailed in \cref{fig:specimen}, with a nominal line spacing $l=\SI{0.4}{mm}$.
Both algorithms are executed on the same machine, a Windows computer with an Intel Core i9-9900K CPU running at \SI{3.60}{\giga\hertz} and using \SI{48}{\giga\byte} of RAM. The global trajectory generation algorithm solves a single large optimization problem and requires \SI{22.7}{s} to produce the print trajectories. 
The proposed swarm-based trajectory generation algorithm solves iteratively a large number of simpler optimization problems and completes the task in \SI{198}{\milli\second}. Thanks to its speed, our approach can be smoothly integrated in existing slicers without affecting user experience. The \SI{99}{\percent} reduction in computation time is only one advantage of the method we propose in this work.
As discussed by the authors in \cite{guidetti2023stress}, the global trajectory generation method returns only one solution for a given part and load case, producing trajectories at a fixed (and quasi constant) distance. This is well suited to printing processes where the deposited line width cannot be changed, but penalizes stress alignment. In fact, the trajectories produced with the global method (also reported in \cref{tab:K}) have $\bar{\beta} = 0.983$ and a variance of normalized distances of \SI{4.4e-4}{}. The trajectory spacing is extremely uniform, and the stress alignment is comparable with the results of \cref{alg:outline} with $K=0.5$. When variable line width printing is possible, the swarm-based trajectory generation method offers more flexibility, as trajectories can be adapted by selecting the correct value of $K$. Allowing a certain amount of variation in the trajectories spacing increases stress alignment (as shown in \cref{tab:K}) and as a consequence improves the mechanical properties of the part. The tuning of $K$ is however limited by the hardware effectiveness in variable width printing. The stress alignment and trajectory spacing metrics we have proposed are a useful tool for choosing the most suitable value of $K$ in practice. Based on these metrics, the ideal $K$ can be found with numerous sampling based methods, such as for example random sampling, grid search, or Bayesian optimization.

\subsection{Printed Parts}

\begin{figure*}[htbp]
\centering
{\scalebox{1}[-1]{\includegraphics[trim={55mm 8mm 55mm 4mm},clip,angle=-90,width=\textwidth]{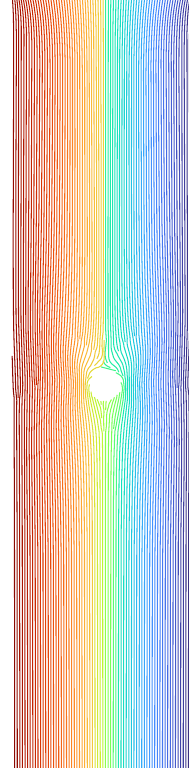}}}
\caption{Print trajectories for the open-hole tensile specimen obtained with swarm-based trajectory generation (\cref{alg:outline}) and $K=5$. Agents move from the left to the right of the figure.}
\label{fig:print_traj_5}
\end{figure*}

\begin{figure}[htbp]
\centering
\includegraphics[]{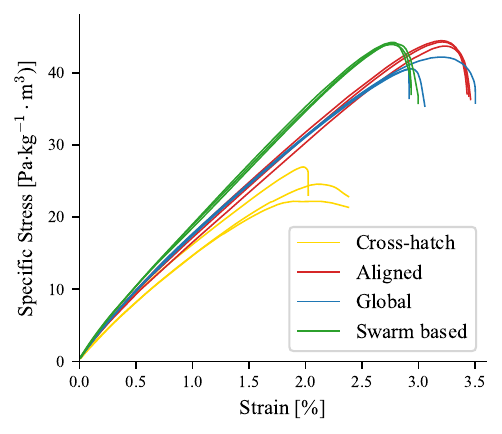}
\caption{Specific stress-strain curve of the twelve tested open-hole tensile specimen}
\label{fig:results}
\end{figure}

We manufacture and test according to \cite{american2023standard} \num{12} open-hole tensile specimen which are printed following different trajectories. The utilized trajectories are generated with:
\begin{enumerate}
    \item A commercial slicer \cite{prusaslicer} set to \emph{cross-hatching} infill. This infill creates a rectilinear grid by printing one layer as parallel lines in one direction, the next layer rotated by \SI{90}{\degree}, etc. The direction of the lines is set to \SI{45}{\degree} with respect to the sides of the specimen. Lines are generated at a distance of \SI{0.4}{mm} and with a \SI{100}{\percent} infill. This case constitutes a baseline for non stress-aligned printing.
    \item A commercial slicer \cite{prusaslicer} set to \emph{aligned rectilinear} infill. All lines in the part are parallel and aligned in the direction of the long side of the specimen. Lines are generated at a distance of \SI{0.4}{mm} and with a \SI{100}{\percent} infill. This case constitutes a baseline for naive stress-aligned printing.
    \item The global trajectory generation method from \cite{guidetti2023stress}. The nominal line distance is set to \SI{0.4}{mm}. This case is a state-of-the-art benchmark.
    \item The proposed swarm-based trajectory generation method with $K=5$ and a nominal line distance of \SI{0.4}{mm}. These print trajectories are shown in \cref{fig:print_traj_5}.
\end{enumerate}
We manufacture three samples for each trajectory generation case. All parts are printed in PLA on a Prusa i3 MK3S machine at a nozzle temperature of \SI{205}{\celsius}, a bed temperature of \SI{60}{\celsius} and feed rate of \SI{3150}{\milli\meter\per\minute}. Samples are tested on a Galdabini Quasar 10 testing machine following \cite{american2023standard}. As different trajectory generation methods produce different line spacing distributions, and as this affects the deposition process, the samples have different densities. To make the results comparable, we divide the force measurements by the density of each sample to produce specific stress, specific modulus and specific strength (which are the density adjusted equivalents of stress, Young's modulus and ultimate tensile strength). The results are reported in \cref{fig:results} and \cref{tab:results}. Cross-hatching, the most common trajectory generation technique in FFF, produces the worst performance, with the lowest specific modulus and strength. The aligned rectilinear approach and the global trajectory generation method have a very similar stress-strain curve, and a comparable specific modulus. The samples manufactured with our swarm-based trajectories outperform all other samples on the entire strain range. They produce a $\sim 10\%$ improvement in specific modulus with respect to the aligned rectilinear and the global methods, while retaining the specific strength of the aligned rectilinear method. Contrarily to previous works \cite{guidetti2023stress,fang2020reinforced}, this improvement was achieved on a conventional planar printer, using a ubiquitous feedstock material. It indicates that any combination of commercially available slicer and printer can benefit from our approach to produce stiffer and stronger parts via optimized trajectory generation. We point out that the high specific strength of the naive aligned rectilinear approach is in all likelihood the consequence of the line spacing homogeneity, which enables high quality deposition. We expect that the constant developments in the field of variable width FFF will reduce the number of deposition defects in the optimized parts, further increasing their strength. Furthermore, the aligned rectilinear approach can be used exclusively in very simple geometries where the load is rectilinear and stress-alignment can be achieved intuitively. Conversely, rigorous trajectory optimization methods can be applied to geometries and load cases of any complexity, as it was demonstrated in \cite{guidetti2023stress}.

\begin{table*}[htbp]
\centering
\caption{Specific modulus and specific strength of the twelve tested open-hole tensile specimen. The values were obtained by averaging the results from samples produced with the same trajectory generation method.}
\renewcommand{\arraystretch}{1.3}
\renewcommand\cellset{\renewcommand\arraystretch{0.8}%
\setlength\extrarowheight{0pt}}
\begin{tabular}[h]{@{}l c c c c @{}}
\toprule
& \makecell[c]{Cross-\\hatching} &
\makecell[c]{Aligned\\rectilinear} &
\makecell[c]{Global method\\from \cite{guidetti2023stress}} &
\makecell[c]{Swarm-based\\method\\($K=5$)} \\

\midrule

\makecell[c]{Specific modulus\\ $[$\SI{}{\pascal\per\kilo\gram\cubic\meter}$]$} & 
\num{17.48} & \num{19.53} & \num{19.91} &  $\mathbf{21.73}$ \\ [1.2em]
\makecell[c]{Specific strength\\ $[$\SI{}{\pascal\per\kilo\gram\cubic\meter}$]$} & 
\num{24.57} & \num{44.14} & \num{41.08} &  $\mathbf{44.08}$ \\
\bottomrule
\end{tabular}
\label{tab:results}
\end{table*}

\section{Conclusion} \label{sec:concl}

In this work, we introduced a novel swarm-based approach to the generation of optimized stress-aligned print trajectories for 3D printing. We evaluated our method and compared with state-of-the-art approaches, both in simulation and in experiments, focusing on the implications of stress-aligned trajectories for mechanical properties and printing quality. Our results demonstrate a significant improvement in computational efficiency with our proposed algorithm, achieving a remarkable $115\times$ increase in computation performance compared to the state-of-the-art method. This computational advantage, coupled with the flexibility of adapting trajectory spacing through hyperparameter tuning, positions our approach as a highly practical solution for seamless integration into existing slicers without compromising user experience.
Furthermore, our method allows for enhanced stress alignment, leading to improved mechanical properties of printed parts. The specific modulus of parts produced using our swarm-based trajectories exhibited a $\sim 10\%$ enhancement compared to existing methods.
Our findings highlight the potential of trajectory optimization in advancing the mechanical performance and efficiency of 3D printing processes in general, and of the prevalent planar FFF of PLA in particular.
In future research, we plan to extend the swarm-based approach to the non-planar slicing problem, and to solve the slicing and trajectory generation tasks in a single optimization. 

\section*{Acknowledgment}
We gratefully acknowledge the support of NematX AG which has provided the equipment to manufacture the samples.

\bibliographystyle{IEEEtran}
\bibliography{ref.bib}

\EOD

\end{document}